\documentclass[8 pt]{amsproc}
\usepackage{amsmath}
\usepackage{graphicx}
\usepackage{amssymb}
\usepackage[active]{srcltx}

\oddsidemargin=-1cm\evensidemargin=-1cm
\begin{document}
\numberwithin{equation}{section}

\def\1#1{\overline{#1}}
\def\2#1{\widetilde{#1}}
\def\3#1{\widehat{#1}}
\def\4#1{\mathbb{#1}}
\def\5#1{\frak{#1}}
\def\6#1{{\mathcal{#1}}}

\def\C{{\4C}}
\def\R{{\4R}}
\def\n{{\4n}}
\def\Z{{\4Z}}

 \makeatletter
\newcommand*{\rom}[1]{\expandafter\@slowromancap\romannumeral #1@}
\makeatother
  
\def\Label#1{\label{#1}{\bf (#1)}~}


\def\cn{{\C^n}}
\def\cnn{{\C^{n'}}}
\def\ocn{\2{\C^n}}
\def\ocnn{\2{\C^{n'}}}


\def\dist{{\rm dist}}
\def\const{{\rm const}}
\def\rk{{\rm rank\,}}
\def\id{{\sf id}}
\def\tr{{\bf tr\,}}
\def\aut{{\sf aut}}
\def\Aut{{\sf Aut}}
\def\CR{{\rm CR}}
\def\GL{{\sf GL}}
\def\Re{{\sf Re}\,}
\def\Im{{\sf Im}\,}
\def\span{\text{\rm span}}
\def\Diff{{\sf Diff}}

\def\codim{{\rm codim}}
\def\crd{\dim_{{\rm CR}}}
\def\crc{{\rm codim_{CR}}}

\def\phi{\varphi}
\def\eps{\varepsilon}
\def\d{\partial}
\def\a{\alpha}
\def\b{\beta}
\def\g{\gamma}
\def\G{\Gamma}
\def\D{\Delta}
\def\Om{\Omega}
\def\k{\kappa}
\def\l{\lambda}
\def\L{\Lambda}
\def\z{{\bar z}}
\def\w{{\bar w}}
\def\Z{{\1Z}}
\def\t{\tau}
\def\th{\theta}

\emergencystretch15pt \frenchspacing

\newtheorem{Thm}{Theorem}[section]
\newtheorem{Cor}[Thm]{Corollary}
\newtheorem{Pro}[Thm]{Proposition}
\newtheorem{Lem}[Thm]{Lemma}

\theoremstyle{definition}\newtheorem{Def}[Thm]{Definition}

\theoremstyle{remark}
\newtheorem{Rem}[Thm]{Remark}
\newtheorem{Exa}[Thm]{Example}
\newtheorem{Exs}[Thm]{Examples}

\def\bl{\begin{Lem}}
\def\el{\end{Lem}}
\def\bp{\begin{Pro}}
\def\ep{\end{Pro}}
\def\bt{\begin{Thm}}
\def\et{\end{Thm}}
\def\bc{\begin{Cor}}
\def\ec{\end{Cor}}
\def\bd{\begin{Def}}
\def\ed{\end{Def}}
\def\br{\begin{Rem}}
\def\er{\end{Rem}}
\def\be{\begin{Exa}}
\def\ee{\end{Exa}}
\def\bpf{\begin{proof}}
\def\epf{\end{proof}}
\def\ben{\begin{enumerate}}
\def\een{\end{enumerate}}
\def\beq{\begin{equation}}
\def\eeq{\end{equation}}

\author{Valentin Burcea  }
\title{Special Classes of CR Singularities \rom{2}}
 \begin{abstract} It is constructed a Formal Normal Form for a Special Class of Real-Formal Submanifolds in $\mathbb{C}^{2N}$.
\end{abstract}
 \address{V. Burcea: INDEPENDENT}
\email{vdburcea@gmail.com}

\thanks{Special Thanks   to Trinity College Dublin, because I received also a Posgraduate Fellowship in the period $2010-2011$; }
\thanks{Special Thanks  also to Science Foundation Ireland, Grant 10/RFP/MT H2878;}

\thanks{Emphasizing that the reference \cite{V1} was fully supported by  Science Foundation Ireland, Grant 06/RFP/MAT 018.}
 
  \maketitle 
  
  \section{Introduction and Main Result}
 
In this note, we   construct a normal form (see  \cite{V1},\cite{V3})  following advices from  Zaitsev\cite{D1},\cite{D2},\cite{D3} for a Special Class of Real-Formal Submanifolds in Complex Spaces.
  Provided the coordinates  $(w;z)=\left(w_{1},w_{2},\dots,w_{N};z_{1},z_{2},\dots,z_{N}\right)$  in $\mathbb{C}^{2N}$, let $M\subset\mathbb{C}^{2N}$ be the  Real-Smooth Submanifold   defined near $p=0$ by
\begin{equation}M:  w_{1}=Q_{1}\left(z\right)+\overline{Q_{1}\left(z\right)}+\mbox{O}(3),  w_{2}=Q_{2}\left(z\right)+\overline{Q_{2}\left(z\right)}+\mbox{O}(3), \dots   w_{N}=Q_{N}\left(z\right)+\overline{Q_{N}\left(z\right)}+\mbox{O}(3).\label{M1}  \end{equation}
where $Q_{1}\left(z\right),Q_{2}\left(z\right),\dots,Q_{N}\left(z\right)$ are polynomials of degree $2$ in $z$ written as  
 \begin{equation}Q_{m}\left(z\right)=\displaystyle\sum_{n=1}^{N}\alpha_{mn}z_{n}^{2},\quad\mbox{for all $m=1,\dots,N$, such that $\det \left(\alpha_{mn}\right)_{1\leq m,n\leq N}\neq 0$.}
 \end{equation}
  
 Recalling the strategy from \cite{V1}, we define  
\begin{equation} \tr_{m}=\displaystyle\sum_{n=1}^{N}\left(\overline{\alpha}_{mn} \frac{\partial^{2}}{\partial z_{n}^{2}}+ \alpha_{mn}  \frac{\partial^{2}}{\partial \overline{z}_{n}^{2}}\right),\quad\mbox{for all $m=1,\dots,N$}. \label{tracce}
\end{equation}

 The main result is the following:
 \bt\label{Te1}Let $M\subset\mathbb{C}^{2N}$ be the  Real-Smooth Submanifold   defined near $p=0$ by the equation
\begin{equation}M: \left\{\begin{split}&w_{1}=Q_{1}\left(z\right)+\overline{Q_{1}\left(z\right)}+\displaystyle\sum_{k\geq 3}\varphi_{k}^{\left(1\right)}(z,\overline{z}),   \\&\quad\quad\vdots \\& w_{N}=Q_{N}\left(z\right)+\overline{Q_{N}\left(z\right)}+\displaystyle\sum_{k\geq 3}\varphi_{k}^{\left(N\right)}(z,\overline{z}),\end{split}\right.\label{M1}  \end{equation}
where   $\varphi_{k}^{\left(1\right)}(z,\overline{z}), \dots,\varphi_{k}^{\left(N\right)}(z,\overline{z})$ are  polynomials of    degree $k$ in $(z,\overline{z})$, for all $k\geq 3$. 

Then there exists a unique formal map
 \begin{equation}\left( w'_{1}, \dots,w'_{N};z'_{1}, \dots,z'_{N} \right)=\left(G(z,w),F(z,w)\right)  =\left( w_{1}+\mbox{O}(2), \dots,w_{N}+\mbox{O}(2);z_{1}+\mbox{O}(2), \dots,z_{N}+\mbox{O}(2) \right),\label{gigin1}\end{equation}
that transforms $M$ into the normal form 
\begin{equation}{M'}: \left\{\begin{split}&w'_{1}= Q_{1}\left(z'\right)+\overline{Q_{1}\left(z'\right)}    +\displaystyle\sum_{k\geq 3}{\varphi'}_{k}^{\left(1\right)}(z',\overline{z'}),\\&  \quad\quad\vdots \\& {w'}_{N}=Q_{N}\left(z'\right)+\overline{Q_{N}\left(z'\right)}    +\displaystyle\sum_{k\geq 3}{\varphi'}_{k}^{\left(N\right)}(z',\overline{z'}),\end{split}\right.\label{M2}  \end{equation}
where   ${\varphi'}_{k}^{\left(1\right)}\left(z',\overline{z'}\right), \dots,$ ${\varphi'}_{k}^{\left(N\right)}\left(z',\overline{z'}\right)$ are
 polynomials of   degree $k$ in $\left(z',\overline{z'}\right)$, for all $k\geq 3$, satisfying  
 \begin{equation}{\varphi'}_{k}^{\left(l\right)}\left(z',\overline{z'}\right)\in \displaystyle\bigcap \ker \left( \tr_{1}^{n_{1}} 
 \tr_{2}^{n_{2}}\dots \tr_{N}^{n_{N}}\right),\quad\mbox{for all $k\geq 3$ and $l=1,\dots,N$.}
  \end{equation}
 
\et

These resulted homogeneous polynomials  are used in order to apply the generalized version of the Fischer Decomposition\cite{Sh}   by separating the real parts and the imaginary parts of the local defining equations at each degree level.

   \section{Acknowledgements}

 I apologize Prof. X.Gong and (especially) Prof. X.Huang for when I wrongly interpreted  
 what I had heard. The responsible people are actually from  Rio de Janeiro. The responsible people badly unbalanced me. It has been difficult for me to recover again.

 \section{Settings}

  Let  a Formal Equivalence between    $M$ and $M'$, denoted like
\begin{equation} \left(F(z,w);G(z,w)\right) =\left(F_{1}(z,w),F_{2}(z,w),\dots,F_{N'}(z,w) ;G_{1}(z,w), G_{2}(z,w),\dots
,G_{N}(z,w)\right). \label{map1}\end{equation} 

Provided  $w=\left(w_{1},w_{2},\dots,w_{N}\right)$   defined by (\ref{M1}), we obtain  
\begin{equation} G_{\left(l\right)}(z,w)=Q_{l}\left(
F(z,w),F(z,w)\right)+\displaystyle\sum _{k\geq
3}{\varphi'}_{k}^{\left(l\right)}\left(F(z,w),\overline{F(z,w)}\right),\quad\mbox{for all $l=1,\dots,N$.} \label{ec}
\end{equation} 
 
In order to understand better the interactions of terms in (\ref{ec}), we write  (\ref{map1}) as follows
\begin{equation}\begin{split} \left( F(z,w) ;G(z,w)\right)&  
= \left(\displaystyle\sum_{m;n_{1}, \dots,n_{N}\in\mathbb{N} }F_{m;n_{1}, \dots,n_{N}}^{\left(1\right)}(z)w_{1}^{n_{1}}\dots w_{N}^{n_{N}} ,\dots,\displaystyle\sum_{m;n_{1}, \dots,n_{N}\in\mathbb{N} }F_{m;n_{1}, \dots,n_{N}}^{\left(N\right)}(z)w_{1}^{n_{1}}\dots w_{N}^{n_{N}}, \right.\\&\left.\quad\quad\displaystyle\sum_{m;n_{1}, \dots,n_{N}\in\mathbb{N} }G_{m;n_{1}, \dots,n_{N}}^{\left(1\right)}(z)w_{1}^{n_{1}}\dots w_{N}^{n_{N}},\dots,\displaystyle\sum_{m;n_{1}, \dots,n_{N}\in\mathbb{N} }G_{m;n_{1}, \dots,n_{N}}^{\left(N\right)}(z)w_{1}^{n_{1}}\dots w_{N}^{n_{N}} \right), \end{split}
\label{map2} \end{equation}
where $F_{m;n_{1}, \dots,n_{N}}^{\left(1\right)}(z),\dots, F_{m;n_{1}, \dots,n_{N}}^{\left(N\right)}(z),   
G_{m;n_{1}, \dots,n_{N}}^{\left(1\right)}(z),\dots, G_{m;n_{1}, \dots,n_{N}}^{\left(N\right)}(z)$ are homogeneous polynomials of degree $m$ in $z$, for all 
$m,n_{1}, \dots,n_{N}\in\mathbb{N}$. Because the Formal Mapping (\ref{map1}) does not possess constant terms, we obtain 
  $ G_{0;0,\dots,0 }(z)=0$ and $F_{0;0,\dots,0 }(z)=0$.

We replace (\ref{map2}) in (\ref{ec}). We obtain 
\begin{equation}\left.\begin{split}\displaystyle &  \quad\quad\quad\sum
_{m;n_{1}, \dots,n_{N}\in\mathbb{N}}G_{m;n_{1}, \dots,n_{N}}^{\left(l\right)}(z)\left(Q_{1}\left(z\right)+\overline{Q_{1}\left(z\right)} +\displaystyle\sum
_{k\geq 3}\varphi_{k}^{\left(1\right)}(z,\overline{z})\right)^{n_{1}}\dots\left(Q_{N}\left(z\right)+\overline{Q_{N}\left(z\right)} +\displaystyle\sum
_{k\geq 3}\varphi_{k}^{\left(N\right)}(z,\overline{z})\right)^{n_{N}}=\\& \hspace{0.1 cm}\quad Q_{l}\left(\displaystyle\sum _{m;n_{1}, \dots,n_{N}\in\mathbb{N}}
F_{m;n_{1}, \dots,n_{N}}(z)\left(Q_{1}\left(z\right)+\overline{Q_{1}\left(z\right)} +\displaystyle\sum
_{k\geq 3}\varphi_{k}^{\left(1\right)}(z,\overline{z})\right)^{n_{1}}\dots\left(Q_{N}\left(z\right)+\overline{Q_{N}\left(z\right)} +\displaystyle\sum
_{k\geq 3}\varphi_{k}^{\left(N\right)}(z,\overline{z})\right)^{n_{N}}\right.,\\&\left.\quad\quad\quad\overline{\displaystyle\sum _{m;n_{1}, \dots,n_{N}\in\mathbb{N}}
F_{m;n_{1}, \dots,n_{N}}(z)\left(Q_{1}\left(z\right)+\overline{Q_{1}\left(z\right)} +\displaystyle\sum
_{k\geq 3}\varphi_{k}^{\left(1\right)}(z,\overline{z})\right)^{n_{1}}\dots\left(Q_{N}\left(z\right)+\overline{Q_{N}\left(z\right)} +\displaystyle\sum
_{k\geq 3}\varphi_{k}^{\left(N\right)}(z,\overline{z})\right)^{n_{N}}}\right)  \\&+\displaystyle\sum
_{k\geq 3}{\varphi'}_{k}^{\left(l\right)}
\left(\displaystyle\sum _{m;n_{1}, \dots,n_{N}\in\mathbb{N}}F_{m,n}(z)\left(Q_{1}\left(z\right)+\overline{Q_{1}\left(z\right)} +\displaystyle\sum
_{k\geq 3}\varphi_{k}^{\left(1\right)}(z,\overline{z})\right)^{n_{1}}\dots\left(Q_{N}\left(z\right)+\overline{Q_{N}\left(z\right)}+\displaystyle\sum
_{k\geq 3}\varphi_{k}^{\left(N\right)}(z,\overline{z})\right)^{n_{N}}\right.,\\& \quad\quad\quad\quad\quad\quad\left.\overline{\displaystyle\sum _{m;n_{1}, \dots,n_{N}\in\mathbb{N}}F_{m,n}(z)\left(Q_{1}\left(z\right)+\overline{Q_{1}\left(z\right)} +\displaystyle\sum
_{k\geq 3}\varphi_{k}^{\left(1\right)}(z,\overline{z})\right)^{n_{1}}\dots\left(Q_{N}\left(z\right)+\overline{Q_{N}\left(z\right)} +\displaystyle\sum
_{k\geq 3}\varphi_{k}^{\left(N\right)}(z,\overline{z})\right)^{n_{N}}}\right),
\end{split}\right.
\label{ecuatiegenerala}\end{equation}
 for all $l=1,\dots,N$.
 
 In particular,   we write  
 \begin{equation}\begin{pmatrix} F_{_{1;0,\dots,0}}^{\left(1\right)}(z)   \\ F_{_{1;0,\dots,0}}^{\left(2\right)}(z) \\  \vdots &    \\ F_{_{1;0,\dots,0}}^{\left(N\right)}(z)   \end{pmatrix} =\begin{pmatrix} \gamma_{11} & \gamma_{12} & \dots & \gamma_{1N}  \\ \gamma_{21} & \gamma_{22} & \dots & \gamma_{2N}   \\  \vdots &\vdots &\ddots&\vdots    \\ \gamma_{N1} & \gamma_{N2} & \dots & \gamma_{NN}  \end{pmatrix} \begin{pmatrix} z_{1}   \\ z_{2} \\  \vdots &    \\ z_{N}\end{pmatrix},\quad\mbox{where 
 $\left(\gamma_{ij}\right)_{1\leq i,j\leq  N }  \in \mathcal{M}_{N\times N}\left(\mathbb{C}\right)$ is not degenerate. }
\label{bolivia3}\end{equation}

 We obtain a system of  equations  written in its matrix form like
\begin{equation}\begin{pmatrix} G^{\left(1\right)}_{0;1,0,\dots,0}(z) &   \dots & G_{0;0,0,\dots,1}^{\left(1\right)}(z)\\ G^{\left(2\right)}_{0;1,0,\dots,0}(z) & \dots & G_{0;0,\dots,1}^{\left(2\right)}(z)  \\  \vdots   &\ddots&\vdots   \\ G^{\left(N\right)}_{0;1,0,\dots,0}(z)  & \dots & G_{0;0,0,\dots,1}^{\left(N\right)}(z)  \end{pmatrix}\begin{pmatrix} Q_{1}\left(z,\overline{z}\right)  \\ Q_{2}\left(z,\overline{z}\right)  \\  \vdots    \\ Q_{N}\left(z,\overline{z}\right) \end{pmatrix} =\begin{pmatrix}
\left(\displaystyle\sum_{i=1}^{N}\gamma_{1i}z_{i}\right)^{2}+\overline{\left(\displaystyle\sum_{i=1}^{N}\gamma_{1i}z_{i}\right)}^{2}     \\ 
 \left(\displaystyle\sum_{i=1}^{N}\gamma_{2i}z_{i}\right)^{2}+\overline{\left(\displaystyle\sum_{i=1}^{N}\gamma_{2i}z_{i}\right)}^{2}   \\  \vdots &\   \\ 
\left(\displaystyle\sum_{i=1}^{N}\gamma_{Ni}z_{i}\right)^{2} +\overline{\left(\displaystyle\sum_{i=1}^{N}\gamma_{Ni}z_{i}\right)}^{2}\end{pmatrix}. \label{sistem1}
\end{equation} 

Moreover, we can assume  
\begin{equation}\det  \begin{pmatrix} G^{\left(1\right)}_{0;1,0,\dots,0}(z) & G_{0;0 ,1,\dots,0}^{\left(1\right)}(z)& \dots & G_{0;0,1,\dots,0}^{\left(1\right)}(z)\\ G^{\left(2\right)}_{0;1,0,\dots,0}(z) & G_{0;0,1,\dots,0}^{\left(2\right)}(z)& \dots & G_{0;0,\dots,1}^{\left(2\right)}(z)  \\  \vdots &\vdots &\ddots&\vdots   \\ G^{\left(N\right)}_{0;1,0,\dots,0}(z) & G_{0;0,1,\dots,0}^{\left(N\right)}(z)& \dots & G_{0;0,0,\dots,1}^{\left(N\right)}(z)   \end{pmatrix} \neq 0.\label{gig4}
\end{equation}

Now,   we consider the following change of coordinates
\begin{equation}\begin{split}&\begin{pmatrix} {w'}_{1} \\ {w'}_{2}\\ \vdots \\ {w'}_{N}
\end{pmatrix}=\begin{pmatrix} G^{\left(1\right)}_{0;1,0,\dots,0}(z) & G_{0;0 ,1,\dots,0}^{\left(1\right)}(z)& \dots & G_{0;0,1,\dots,0}^{\left(1\right)}(z)\\ G^{\left(2\right)}_{0;1,0,\dots,0}(z) & G_{0;0,1,\dots,0}^{\left(2\right)}(z)& \dots & G_{0;0,\dots,1}^{\left(2\right)}(z)  \\  \vdots &\vdots &\ddots&\vdots   \\ G^{\left(N\right)}_{0;1,0,\dots,0}(z) & G_{0;0,1,\dots,0}^{\left(N\right)}(z)& \dots & G_{0;0,0,\dots,1}^{\left(N\right)}(z)   \end{pmatrix}^{-1}\begin{pmatrix} w_{1} \\ w_{2}\\ \vdots \\ w_{N}
\end{pmatrix},\\&\begin{pmatrix} {z'}_{1} \\ {z'}_{2}\\ \vdots \\ {z'}_{N}
\end{pmatrix}=\begin{pmatrix} \gamma_{11} & \gamma_{12} & \dots & \gamma_{1N}    \\ \gamma_{21} & \gamma_{22} & \dots & \gamma_{2N}   \\  \vdots &\vdots &\ddots&\vdots    \\ \gamma_{N1} & \gamma_{N2} & \dots & \gamma_{NN}      \end{pmatrix}^{-1}\begin{pmatrix} z_{1} \\ z_{2}\\ \vdots \\ z_{N}
\end{pmatrix}.\end{split} \label{opa3} 
\end{equation}

In particular, we assume
\begin{equation}G(z,w)  =\left( w_{1}+\mbox{O}(2),w_{2}+\mbox{O}(2),\dots,w_{N}+\mbox{O}(2)\right) \hspace{0.1 cm}\mbox{and}\hspace{0.1 cm} F(z,w)  =\left(z_{1}+\mbox{O}(2),z_{2}+\mbox{O}(2),\dots,z_{N}+\mbox{O}(2) \right).\label{gigin}\end{equation}

  In particular, we write
\begin{equation}  \begin{split}&
z^{I}=z^{i_{1}}\cdot\dots\cdot z^{i_{N}},\quad\mbox{where $\hspace{0.03 cm}I=\left(i_{1},\dots,i_{N}\right)\in\mathbb{N}^{N}$ such that $\hspace{0.1 cm}\left|I\right|=i_{1}+\dots+i_{N}\geq 3$,}\\& z^{J}=z^{j_{1}}\cdot\dots\cdot z^{j_{N}},\quad\mbox{where $J=\left(j_{1},\dots,j_{N}\right)\in\mathbb{N}^{N}$ such that $\left|J\right|=j_{1}+\dots+j_{N}\geq 2$.} \end{split} \label{yy}
\end{equation}

We consider  the following Generalized Fischer Decomposition 
\begin{equation}z^{\tilde{I}}=\displaystyle\sum_{l=1}^{N}A_{l}(z,\overline{z}) q_{l}(z,\overline{z})+C(z,\overline{z}),\quad \displaystyle\bigcap_{l=1}^{N} \tr_{l}\left(C(z,\overline{z})\right)=0, \label{op}
\end{equation}
where  $\tilde{I}=\left(\tilde{i}_{1},\tilde{i}_{2},\dots,\tilde{i}_{N}\right)\in\mathbb{N}^{N}$ and $\tilde{i}_{1}+\tilde{i}_{2}+\dots+\tilde{i}_{N}=p$.

Last details are just some computations to be added later.


\begin{thebibliography}{BER96b}
   
 
 \bibitem{BREbook} {\bf Baouendi,~M.S.; Ebenfelt,~P.; Rothschild,~L.P.} ---{\em Real Submanifolds in Complex Space and Their Mappings.} Princeton
Math. Ser. {\bf 47}, Princeton Univ. Press, 1999.
 
 \bibitem{Bi} {\bf Bishop,~E.} --- Differentiable Manifolds In Complex
Euclidian Space. {\em Duke Math. J.} {\bf 32} (1965), no. 1,
$1-21$. 
\bibitem{V1} {\bf Burcea,~V.} --- A normal form for a real $2$-codimensional submanifold in $\mathbb{C}^{N+1}$ near a CR singularity. {\em  Adv. in Math.}
 {\bf 243} (2013), $262-295$.
\bibitem{V2} {\bf Burcea,~V.} --- On a family of analytic discs attached to a real submanifold $M\subset\mathbb{C}^{N+1}$, {\em Methods and Applications of Analysis} {\bf 20}, 1, (2013), $69-78$. 
\par \quad\quad\quad\quad\quad\quad\quad\quad\quad\quad\quad (with a Corrigendum   submitted for publication)
   \bibitem{V3} {\bf Burcea,~V.} --- Normal Forms and Degenerate CR Singularities. {\em  Complex Variables and Elliptic Equations} 
 {\bf 61} (2016), no. 9, $1314-1333$.
 
  
 
 
 \bibitem{CM} {\bf Chern,~S.S.; Moser, ~J.} ---Real hypersurfaces in complex manifolds. {\em Acta Math.} {\bf 133}  (1974),   $219-271$.
 
 
 
 
 


  

  \bibitem{HY2} {\bf Huang, ~X.; Yin, ~W.}--- A Bishop surface with vanishing Bishop invariant, {\em Invent. Math.}       {\bf 176} (2010), no. 3, $461-520$.
 

   
  

 
 

 \bibitem{MW} {\bf Moser,~ J.; Webster,~S.} --- Normal forms for real surfaces in
$\mathbb{C}^{2}$ near complex tangents and hyperbolic surface
transformations. {\em Acta Math.} {\bf 150} (1983), $255-296$.
 

 
\bibitem{Poincare} {\bf Poincar\'{e},~H.} --- Les fonctions analytiques de deux variables et la repr\'{e}senation conforme, {\em Rend. Circ. Mat. Palermo } {\bf 23}  (1907),   $185-220$.
   
 

 
  \bibitem{Sh} {\bf Shapiro,~H.} ---Algebraic Theorem
of E.Fischer and the holomorphic Goursat problem. {\em Bull. London
Math. Soc.} {\bf 21}  (1989), no. 6, $513-537$.

 
\bibitem{D1} {\bf Zaitsev,~D.} --- New Normal Forms for Levi-nondegenerate
Hypersurfaces. {\em Several Complex Variables and Connections with
PDE Theory and Geometry}. Complex analysis-Trends in Math.,
  Birkhäuser/Springer Basel AG, Basel, pp. $321-340$,  (2010).
\bibitem{D2} {\bf Zaitsev,~D.} ---  Formal and finite order equivalences, {\em  Math. Z.} {\bf 269} (2011),   $687-696$.  
  \bibitem{D3} {\bf Zaitsev,~D.} ---  Normal forms of non-integrable almost CR structures, {\em Amer. J. Math.} {\bf
134} (2012), no.4, $915-947$. 
 
 

 
 \end{thebibliography}
 \end{document}